\documentclass[11pt]{article}

\usepackage{amsthm}
\usepackage{amsmath}
\usepackage{amssymb}

\usepackage{graphicx}

\newtheorem*{lemma}{Lemma}
\newtheorem*{theorem}{Theorem}
\newtheorem{ntheorem}{Theorem}
\newtheorem*{definition}{Definition}

\newtheorem*{proposition}{Proposition}

\newtheorem{nproposition}{Proposition}

\newcommand{\x}{{\boldsymbol{x}}}

\newcommand{\R}{\mathbb{R}}
\newcommand{\RR}{\mathbb{R}}
\newcommand{\TP}{\mathbb{TP}}

\newcommand{\w}{\omega}

\newcommand{\F}{{\cal F}}

\newcommand{\+}{\oplus}

\renewcommand{\x}{\odot}
\newcommand{\tconv}{\text{tconv}}
\newcommand{\cl}{\text{cl}}

\newcommand{\B}{\widetilde{\cal B}}

\begin{document}
\title{Subdominant matroid ultrametrics}
\author{Federico Ardila}
\date{}
\maketitle

\begin{abstract}

Given a matroid $M$ on the ground set $E$, the Bergman fan
$\B(M)$, or space of $M$-ultrametrics, is a polyhedral complex in
$\R^E$ which arises in several different areas, such as tropical 
algebraic geometry, dynamical systems, and phylogenetics. Motivated
by the phylogenetic situation, we study the following problem:
Given a point $\w$ in $\R^E$, we wish to find an $M$-ultrametric which 
is closest to it in the $\ell_{\infty}$-metric.

The solution to this problem follows easily from the existence of 
the subdominant $M$-ultrametric: a componentwise maximum $M$-ultrametric 
which is componentwise smaller than $\w$. A procedure for computing it is
given, which brings together the points of view of matroid theory
and tropical geometry. 

When the matroid in question is the graphical matroid of the
complete graph $K_n$, the Bergman fan $\B(K_n)$ parameterizes the
equidistant phylogenetic trees with $n$ leaves. In this case, our 
results provide a conceptual explanation for Chepoi and Fichet's 
method for computing the tree that most closely matches measured 
data.

\end{abstract}

\section{Introduction}

Given a matroid $M$ on the ground set $E$, the Bergman fan
$\B(M)$, or space of $M$-ultrametrics, is a polyhedral complex in
$\R^E$ which arises in several different areas, such as tropical 
algebraic geometry \cite{Sturmfels1}, dynamical systems 
\cite{Einsiedler} and phylogenetics \cite{Ardila}. It has been 
described topologically and combinatorially \cite{Ardila}.
Motivated by the phylogenetic situation, we study the following 
problem: Given a point $\w$ in $\R^E$, we wish to find an 
$M$-ultrametric which is closest to it in the $\ell_{\infty}$-metric.

In Section \ref{section:ultrametrics} we define the Bergman fan
$\B(M)$, as well as the notion of an $M$-ultrametric. We offer several
characterizations of them. When $M(K_n)$ is the graphical matroid
of the complete graph $K_n$, $M(K_n)$-ultrametrics are precisely
ultrametrics in the usual sense.

In Section \ref{section:subdominant} we show that, in fact, there
is a componentwise maximum $M$-ultrametric which is componentwise
smaller than $\w$. We call it the \emph{subdominant
$M$-ultrametric} of $\w$, and denote it $\w^M$. A simple
translation of it is a closest $M$-ultrametric to $\w$. A
procedure for computing subdominant $M$-ultrametrics is given,
similar in spirit to Tarjan and Kozen's blue and red rules
\cite{Kozen, Tarjan} for computing the bases of minimum weight of
a matroid.

In Section \ref{section:tropical} we prove that the Bergman fan is
a tropical polytope in the sense of \cite{Develin}, and that the
subdominant $M$-ultrametric $\w^M$ is precisely the tropical
projection of $\w$ onto $\B(M)$.

In Section \ref{section:trees}, we discuss a special case of
particular importance: the Bergman fan of $M(K_n)$, the graphical
matroid of the complete graph $K_n$. As shown by Ardila and
Klivans \cite{Ardila}, the Bergman fan $\B(K_n)$ can be regarded
as a space of phylogenetic trees, and we are interested in finding
an (equidistant) phylogenetic tree that most closely matches
measured data. In this case, our results provide a conceptual 
explanation for the method developed by Chepoi and Fichet in 
\cite{Chepoi} to compute the tree that most closely matches 
measured data.

Throughout this paper, familiarity with the fundamental notions of
matroid theory will be very useful; we refer the reader to
\cite[Ch. 1,2]{Oxley} for the basic definitions.

\section{The Bergman fan and matroid ultrametrics} \label{section:ultrametrics}

Let $M$ be a matroid of rank $r$ on the ground set $E$, and let
$\w \in \R^E$. Regard $\w$ as a weight function on $M$, so that
the weight of a basis $B=\{b_1, \ldots, b_r\}$ of $M$ is given by
$\w_B= \w_{b_1} +\w_{b_2} + \cdots + \w_{b_r}$.

Let $M_{\w}$ be the collection of bases of $M$ having minimum
$\w$-weight. This collection is itself the set of bases of a
matroid; for more information, see for example \cite{Ardila}.

\begin{definition}
The \emph{Bergman fan} of a matroid $M$ with ground set $E$ is:
\[
\B(M) := \{\,  \w \in \R^E \,\, : \,\, \textrm{every element of
$E$ is in a $\w$-minimum basis}\}.
\]
An \emph{$M$-ultrametric} is a vector in $\B(M)$.
\end{definition}

\begin{figure}[h]
\centering
\includegraphics[height=3.5cm]{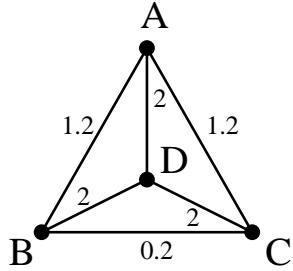}
\caption{An $M(K_4)$-ultrametric.} \label{fig:ultrametric}
\end{figure}

Our ongoing example throughout the paper will be the Bergman fan
of the matroid $M(K_4)$, the graphical matroid of the complete
graph $K_4$. An $M(K_4)$-ultrametric is an assignment of weights
to the edges of $K_4$ such that any edge of $K_4$ is in a spanning
tree of minimum weight. An example of an $M(K_4)$-ultrametric is
given in Figure \ref{fig:ultrametric}. The $\w$-minimum spanning
trees are those consisting of an edge of weight $2$, an edge of
weight $1.2$, and an edge of weight $0.2$. Every edge of the graph
is in at least one such tree.

\medskip

We start by reviewing some useful facts about matroids. Given a
basis $B$ and an element $x \notin B$, there is a unique circuit
which is contained in $B \cup x$ (and must contain $x$). It is
called the \emph{fundamental circuit of $B$ and $x$}.

Given a basis $B$ and an element $y \in B$, there is a unique
cocircuit which is disjoint with $B - y$ (and must contain $y$).
It is called the \emph{fundamental cocircuit of $B$ and $y$}. It
is equal to the fundamental circuit of $E-B$ and $y$ in the dual
matroid $M^*$.

\begin{nproposition}\cite[Lemma 7.3.1]{Bjorner}\label{fundamental}
Let $B$ be a basis of $M$, and let $x \notin B$ and $y \in B$. The
following are equivalent.

\begin{enumerate}
\item[(i)] $B \cup x - y$ is a basis.

\item[(ii)] $y$ is in the fundamental circuit of $B$ and $x$.

\item[(iii)] $x$ is in the fundamental cocircuit of $B$ and $y$.
\end{enumerate}
\end{nproposition}

We will also need the following lemma:

\begin{lemma}\cite[Prop. 2.1.11]{Oxley}\label{circuit-cocircuit}
A circuit and a cocircuit cannot intersect in exactly one element.
\end{lemma}

We now give two additional characterizations of the Bergman fan,
which will be central to our analysis.

\begin{nproposition}\label{Bergman}
Let $M$ be a matroid with ground set $E$, and let $\w \in \R^E$.
The following are equivalent:

\begin{enumerate}

\item[(i)] $\w$ is an $M$-ultrametric.

\item[(ii)] No circuit has a unique $\w$-maximum element.

\item[(iii)] Every element of $E$ is $\w$-minimum in some
cocircuit.
\end{enumerate}

\end{nproposition}

\begin{proof}

$(ii) \Rightarrow (i)$: Let $e \in E$; we want to show that it is
in an $\w$-minimal basis. Consider any $\w$-minimal basis $B$. If
$e \in B$, there is nothing to prove. Otherwise, the fundamental
circuit of $B$ and $e$ has at least two $w$-maximum elements; let
$f \neq e$ be one of them. By Proposition \ref{fundamental}, $B
\cup e - f$ is a basis; since $\w_f \geq \w_e$, it is actually the
$\w$-minimum basis that we need.

$(iii) \Rightarrow (ii)$: Assume that a circuit $C$ achieves its
$\w$-maximum only at $e$. Say $C^*$ is a cocircuit where $e$ is
$\w$-minimum. By the Lemma, we can find an element $f \neq e$ in
$C \cap C^*$; then $\w_e > \w_f$ because $e$ is the unique
$\w$-maximum in $C$, and $\w_e \leq \w_f$ because $e$ is
$\w$-minimum in $C^*$.

$(i) \Rightarrow (iii)$: Let $e \in E$; we want to show that it is
$\w$-minimum in some cocircuit. Let $B$ be an $\w$-minimum basis
containing $e$, and let $f$ be the $\w$-minimum element in the
fundamental cocircuit $C^*$ of $B$ and $e$. Then $B \cup f - e$ is
a basis, and its weight is at least $\w_B$. Therefore $\w_f \geq
\w_e$, so $e$ is $\w$-minimum in $C^*$ also.

\end{proof}

Let us check Proposition \ref{Bergman} for the
$M(K_4)$-ultrametric of Figure \ref{fig:ultrametric}. Statement
$(ii)$ is clear: in each cycle of $K_4$, the two largest weights
are equal. To check statement $(iii)$, recall that the cocircuits
of $M(K_4)$ are the cuts of $K_4$. Denote by $S-S'$ the cut that
separates the vertices in $S$ from the vertices in $S'$. Then the
edges of weight $2$ are minimum in the cut $ABC-D$, the edges of
weight $1.2$ are minimum in the cut $A-BCD$, and the edge of
weight $0.2$ is minimum in the cut $AB-CD$.

\section{Subdominant $M$-ultrametrics.}
\label{section:subdominant}

For $\w,\w' \in \R^E$, say that $\w \leq \w'$ if $\w_e \leq \w'_e$
for each $e \in E$.

\begin{proposition} \label{subdominant}
Let $M$ be a matroid on $E$, and let $\w \in \R^E$. There exists a
unique maximum $M$-ultrametric which is less than or equal to
$\w$. We call it the \emph{subdominant $M$-ultrametric of $\w$},
and denote it $\w^M$.
\end{proposition}

\begin{proof}
Let $S = \{\w' \in \B(M) \,\, : \,\, \w' \leq \w\}$, and let $\w^M
= \sup S$, where the $\sup$ is taken componentwise. We claim that
$\w^M \in \B(M)$.

Proceed by contradiction; assume that $e$ is the unique
$\w^M$-maximum of the circuit $C$. Let $\epsilon$ be such that
$\w^M_e - \epsilon > \w^M_f$ for all $f \in C-e$. We can find a
$\w' \in S$ such that $\w'_e > \w^M_e - \epsilon$. But then $\w'_e
> \w^M_f \geq \w'_f$ for all $f \in C-e$, so $e$ is the unique
$\w'$-maximum of $C$. This contradicts the assumption that $\w'
\in \B(M)$.
\end{proof}

We now turn to the problem of constructing the subdominant
$M$-ultra\-metric of a vector $\w$. Take a vector $\w$ which is
not an $M$-ultrametric. This means that conditions $(ii)$ and
$(iii)$ of Proposition \ref{Bergman} are not satisfied.

We must force each element to be $\w$-minimum in some cocircuit,
and not allow it to be the unique $\w$-maximum element of a
circuit. At the same time, we must demand that it keeps its weight
as high as possible. Let us issue the following two rules, similar
in spirit to Tarjan's \emph{blue and red rules} for constructing
the minimum weight spanning trees of a graph. \cite{Kozen, Tarjan}

\begin{flushleft}
\begin{tabular}{p{4.6in}}

\noindent {\bf Blue Rule:} Suppose an element $e$ is not
$\w$-minimum in any cocircuit. Look at the minimum weight in each
cocircuit containing $e$. (They are all less than $\w_e$.) Make
$\w_e$ equal to the largest such weight.

\medskip

\noindent {\bf Red Rule:} Suppose an element $e$ is the unique
$\w$-maximum element of a circuit. Look at the maximum weight in
$C-e$ for each circuit $C$ containing $e$. (At least one of them
is less than $\w_e$.) Make $\w_e$ equal to the smallest such weight. \\

\end{tabular}
\end{flushleft}

\begin{figure}[h]
\centering
\includegraphics[height=3cm, width=13cm]{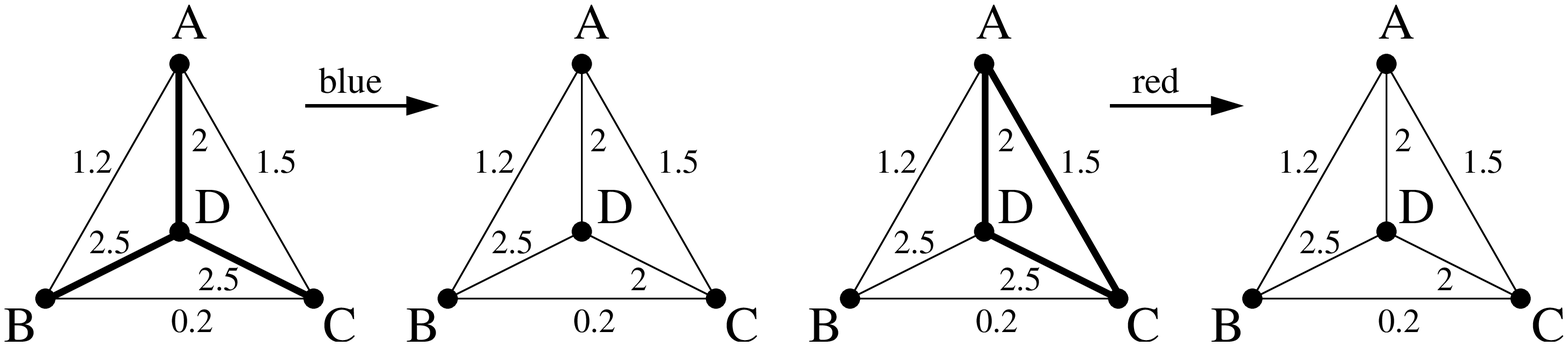}
\caption{Applying the blue rule or the red rule to edge $CD$.}
\label{fig:blueandred}
\end{figure}

Figure \ref{fig:blueandred} shows an example of an assignment $\w$
of weights to the edges of $K_4$ which is not an ultrametric, and
the result of applying the blue rule or the red rule to edge $CD$.
When the blue rule is applied, the edge $CD$ inherits its new
weight from the cut $ABC-D$, as shown. When the red rule is
applied, the edge $CD$ inherits its new weight from either one of
the cycles $ACD$ (as shown) or $ABCD$. Notice that, surprisingly,
the blue rule and the red rule give the same result. The reason
for this will soon be explained.

In principle, it is not yet clear how to apply these rules, or
what they will do. We do not know in what order we should apply
them, and different orders would seem likely to give different
results. It is not even obvious that applying these two local
rules successively will accomplish the global goal of turning $\w$
into an $M$-ultrametric.

Fortunately, as the following theorem shows, the situation is much
simpler than expected. It turns out that, if each element of $E$
is willing to do its part by complying with the two rules, the
global goal of constructing the subdominant $M$-ultrametric will
inevitably be achieved.

\begin{ntheorem}\label{main}
Let $M$ be a matroid on the ground set $E$, let $\w \in \R^E$. For
each element $e \in E$, there are two possibilities:

\begin{enumerate}

\item Both the blue rule and the red rule apply to $e$, and they
both change its weight from $\w_e$ to $\w^M_e$, or

\item Neither the blue rule nor the red rule apply to $e$, and
$\w_e = \w^M_e$.
\end{enumerate}

Consequently, if we apply the blue rule or the red rule to each
element of $E$ in any order, we obtain the subdominant
$M$-ultrametric $\w^M$.
\end{ntheorem}

\begin{figure}[h]
\centering
\includegraphics[height=3cm, width=13cm]{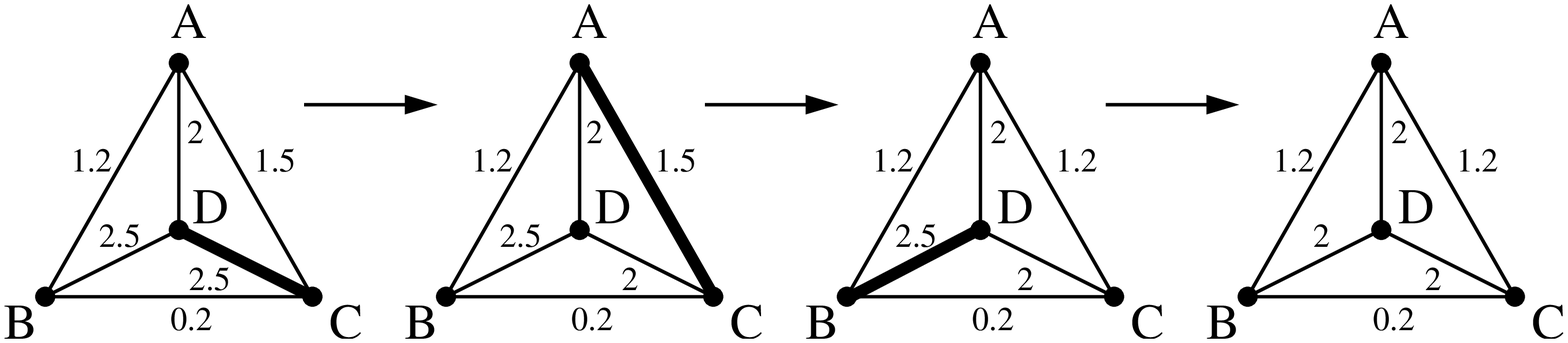}
\caption{Constructing the subdominant $M(K_4)$-ultrametric of
$\w$.} \label{fig:bluered}
\end{figure}

Before proving Theorem \ref{main}, let us illustrate it with an
example: the construction of the subdominant ultrametric of the
weight vector $\w$ we considered in Figure \ref{fig:blueandred}.
Figure \ref{fig:bluered} shows the result of applying, at each
step, either the blue rule or the red rule to the highlighted
edge. The reader is invited to check that, at each step of the
process, the blue rule and the red rule give the same result.

After three steps, we reach an ultrametric. This means that we are
done: we have in fact reached $\w^{M(K_4)}$, the subdominant
$M(K_4)$-ultrametric of $\w$. Moreover, we
could have applied \emph{either the blue rule or the red rule} to
the edges of $K_4$ \emph{in any order}, and we would have obtained
exactly the same result: $\w^{M(K_4)}$.

\begin{proof}[Proof of Theorem \ref{main}]
Let $\w'_e$ and $\w''_e$ be the weights assigned to $e$ by the
blue rule and the red rule, respectively. If the blue rule (or the
red rule) does not apply to $e$, set $\w'_e = \w_e$ (or
$\w''_e=\w_e$). We proceed in several steps.

\medskip

\noindent 1. $\w \geq \w' \geq \w^M$ and $\w \geq \w'' \geq \w^M$.

The inequalities $\w_e \geq \w'_e \geq \w^M_e$ and $\w_e \geq
\w''_e \geq \w^M_e$ are easy to see: When one of the rules
applies, it decrees a decrease in the weight of $e$ that is
clearly necessary for any $M$-ultrametric which is less than $\w$.

\medskip

\noindent 2. $\w'$ is an ultrametric.

We use Proposition \ref{Bergman}$(iii)$. Let $e \in E$. By
definition, $\w'_e$ is equal to the $\w$-minimum weight of some
cocircuit $C^*$ containing $e$. (This is true even if the blue
rule does not apply.) We claim that it is also the $\w'$-minimum
weight of $C^*$. Let $f \in C^*$. Then $\w'_f$ is the largest
$\w$-minimum weight of a cocircuit containing $f$; this is greater
than or equal to the $\w$-minimum weight of $C^*$, which is
$\w'_e$.

\medskip

\noindent 3. $\w' = \w^M$.

We already knew that $\w' \geq \w^M$. Now, since $\w'$ is an
ultrametric less than or equal to $\w$, it must be less than or
equal to $\w^M$ by Proposition \ref{subdominant}.

\medskip

\noindent 4. $\w'' = \w^M$.

We already knew that $\w'' \geq \w^M$. Now assume that $\w''_e > t
> \w'_e$. Let $S = \{f \in E \, | \, \w_f < t\}$. We have
$\w_e \geq \w''_e>t$, so $e$ is not in $S$. Since $\w''_e > t$, no
circuit containing $e$ is in $S \cup e$. Therefore $e \notin
\cl(S)$, and there is a hyperplane $H$ containing $\cl(S)$ with $e
\notin H$. The complement of $H$ is a cocircuit $C^*$ which
contains $e$, and does not contain any element of $S$; its
$\w$-minimum weight is at least $t$. Therefore $\w'_e \geq t$, a
contradiction.
\end{proof}

\medskip

If we have an $\w$-minimum basis of $M$, which is easy to construct
using Tarjan and Kozen's blue and red rules \cite{Tarjan, Kozen}, the 
task of constructing $\w^M$ becomes even simpler, as the following
result shows.

\begin{proposition}\label{basis-pi}
Let $B$ be any $\w$-minimum basis.
\begin{enumerate}
\item[(i)] For each $e \in B$, $\w^M_e = \w_e$.

\item[(ii)] For each $e \notin B$, let $C_0$ be the fundamental
circuit of $B$ and $e$. Then $\w^M_e$ is equal to the $\w$-maximum 
weight in $C_0-e$.
%
\end{enumerate}
\end{proposition}

\begin{proof}
It follows from Proposition \ref{fundamental} that any element $e
\in B$ is $\w$-minimum in the fundamental cocircuit of $B$ and
$e$. Therefore the blue rule does not apply to $e$, and \emph{(i)}
follows from Theorem \ref{main}.

Now we show that, for $e \notin B$, $\min_{e \in C} \max_{j \in
C-e} \, \w_j = \max_{j \in C_0-e} \, \w_j.$ Say the maximum in the
right hand side is achieved at $\w_m$. First we claim that, for
any circuit $C$ containing $e$, there is a basis $B - m \cup f$
with $f \in C-e$. In fact, $B - m \cup C$ has full rank, since it
contains the basis $B - m \cup e$. But $e$ is dependent on $C-e$,
and hence on $(B-m) \cup (C-e)$ also. Therefore $(B-m) \cup (C-e)$
has full rank as well, and the claim follows. Now, since $B - m
\cup f$ is a basis, $\w_m \leq \w_f \leq \max_{j \in C-e} \,
\w_j$. The desired equality follows.

Now \emph{(ii)} follows easily. If the red rule applies to $e$, it
will change $\w_e$ to $\w^M_e = \w_m$. If it does not apply, then
$e$ is not the unique $\w$-max of any circuit, and $\w^M_e = \w_e
= \w_m$.
\end{proof}

\section{Tropical projection} \label{section:tropical}

In this section we study the Bergman complex of a matroid from the
point of view of tropical geometry. We start by reviewing the
basic definitions; for more information, we refer the reader to
\cite{Develin}.

The \emph{tropical semiring} $(\R \cup \{\infty\}, \+, \x)$ is the
set of real numbers augmented by infinity, together with the
operations of \emph{tropical addition and multiplication}, which
are defined by $x \+ y = \min(x,y)$ and $x \x y = x+y$.

Throughout this section, the symbol $\RR$ will denote the tropical
semiring, and not the ring of real numbers. The additive and
multiplicative units of $\RR$ are $\infty$ and $0$, respectively.

The $n$-dimensional space $\RR^n$ is a semimodule over the
tropical semiring, with addition $x \+ y$ and scalar
multiplication $c \x x$ given componentwise. Tropical scalar
multiplication by $c$ is equivalent to translation by the vector
$(c,c,\ldots,c)$. The \emph{$(n-1)$-dimensional tropical
projective space} is $\TP^{n-1} = \RR^n / (x \sim (c \x x))$.

A subset $S$ of $\RR^n$ is \emph{tropically convex} if $(a \x x)
\+ (b \x y) \in S$ for any $x,y \in S$ and \emph{any} $a,b \in
\R$. Tropically convex sets are invariant under tropical scalar
multiplication; \emph{i.e.}, translation by $(c,c,\ldots, c)$. We
will therefore identify them with their image in $\TP^{n-1}$. A
\emph{tropical polytope} is the set
\[
\tconv(V) = \{(a_1 \x v_1) \+ \cdots \+ (a_r \x v_r)\,\, : \,\,
a_1, \ldots, a_r \in \RR \}
\]
of all linear combinations of a finite set $V = \{v_1, \ldots,
v_r\} \subseteq \TP^{n-1}$.

Our claim is that, for any matroid $M$ on the set $[n]$, $-\B(M)$
is a tropical polytope in $\TP^{n-1}$. For each flat $F$ of $M$,
let $v_F \in \TP^{n-1}$ be the vector whose $i$-th coordinate is
$\infty$ if $i \in F$, and $0$ otherwise. Recall that the
\emph{hyperplanes} of a matroid are its maximal proper flats.

\begin{proposition}\label{bergmanistropical}
For any matroid $M$ on $[n]$, $-\B(M)$ is a tropical polytope in
$\TP^{n-1}$. Its set of vertices is
\[
V_M = \{v_H \,\, : \,\, H \text{ is a hyperplane of } M\}.
\]
\end{proposition}

\begin{proof}
We start by reviewing the description of the Bergman fan obtained
in \cite{Ardila}. Given a flag of subsets $ \F=\{\emptyset=:F_0
\subset F_1 \subset \cdots \subset F_{k} \subset F_{k+1}:=E\}, $
the \emph{weight class} of $\F$ is the set of $\w \in \R^n$ for
which $\omega$ is constant on each set $F_i - F_{i-1}$, and has
$\omega|_{F_i-F_{i-1}} < \omega|_{F_{i+1}-F_i}$. For example, one
of the weight classes in $\R^5$ is the set of vectors $\w$ such
that $\w_1=\w_4 < \w_2 < \w_3 = \w_5$. It corresponds to the flag
$\{\emptyset \subset \{1,4\} \subset \{1,2,4\} \subset
\{1,2,3,4,5\}\}$.

The disjoint union of all weight classes is $\R^n$. The Bergman
fan $\B(M)$ is also a union of disjoint weight classes: the weight
class of $\F$ is in $\B(M)$ if and only if $\F$ is a flag of flats
of $M$.

Now we will show that
\[
-\B(M) = \tconv \{v_F \,\, : \,\, F \text{ is a proper flat of }
M\}.
\]
The left hand side is contained in the right hand side because if
$\F$ is a flag of flats, the negative of a vector $w$ in the
weight class of $\F$ can be obtained as
\[
-\w = (-\w|_{F_1} \x v_{F_0}) \+  (-\w|_{F_2 - F_1} \x v_{F_1}) \+
\cdots \+ (-\w|_{F_r - F_{r-1}} \x v_{F_{r-1}}).
\]
To see that the right hand side is contained in $-\B(M)$, since
$v_F \in -\B(M)$ for any flat $F$, it suffices to show that
$-\B(M)$ is tropically convex.

A consequence of the previous description of $\B(M)$ is the
following. A vector $x \in \RR^n$ is in $-\B(M)$ if and only if,
for any $r \in \RR$, the set $(x;r) = \{i \in [n] : x_i \geq r\}$
is a flat of $M$. Now take $x,y \in -\B(M)$ and $a,b \in \RR$, and
let $z = (a \x x) \+ (b \x y)$, so $z_i = \min(a+x_i, b+y_i)$.
Then, for any $r \in \RR$, we have that $(z;r) = (x; r-a) \cap (y;
r-b)$. This is a flat in $M$, because both $(x; r-a)$ and $(y;
r-b)$ are flats. Thus $z \in -\B(M)$.

\medskip

Finally, we prove the claim about the vertices of $-\B(M)$. If $F
= F_1 \cap \cdots \cap F_k$ is an intersection of larger flats,
then
\[
v_F = (0 \x v_{F_1}) \+ \cdots \+ (0 \x v_{F_k})
\]
so $v_F$ is not a vertex. Conversely, suppose that we have an
equation
\[
v_F = (a_1 \x v_{F_1}) \+ \cdots \+ (a_k \x v_{F_k}).
\]
We can assume that $a_i \neq \infty$ for all $i$. For each $f \in
F$, $(v_F)_f = \infty$. Thus for all $i$ we have $(v_{F_i})_f =
\infty$; \emph{i.e.}, $f \in F_i$. For each $\bar{f} \notin F$,
$(v_F)_{\bar{f}} = 0$. Thus for some $i$ we have
$(v_{F_i})_{\bar{f}} \neq \infty$; \emph{i.e.}, $\bar{f} \notin
F_i$. We conclude that $F = F_1 \cap \cdots \cap F_k$.
\end{proof}

For any tropical polytope $P$ in $\TP^{n-1}$, Develin and
Sturmfels \cite{Develin} gave an explicit construction of a
\emph{nearest point map} $\pi_P: \TP^{n-1} \rightarrow P$, which
maps every point in tropical projective space to a point in $P$
which is closest to it in the $\ell_{\infty}$-metric:
\[
||\,x-y||_{\infty} = \max_{1 \leq i,j \leq n} |(x_i-x_j) -
(y_i-y_j)|.
\]

This nearest point map is given as follows. Let $x \in \TP^{n-1}$.
For each vertex $v$ of $P$, we need to compute $\lambda_v$: the
minimum $\lambda$ such that $(\lambda \x v) \+ x = x$. Then
\[
\pi_P(x) = \bigoplus_{v \,\, \textrm{vertex of} \,\, P}
\left(\lambda_v \x v \right).
\]

\begin{nproposition}\label{projection}
The tropical projection $\pi_{-\B(M)}$ maps each vector $\w$ to
its subdominant $M$-ultrametric:
\[
\pi_{-\B(M)}(-\w) = -\w^M.
\]
\end{nproposition}

\begin{proof}

In our case, it is easy to see that $\lambda_{v_H} = \max_{f
\notin H} \w_f$. Therefore
\begin{eqnarray*}
\pi(-\w)_e &=& \min_{\{H \, : \, e \notin H\}} \,\, \max_{f \notin
H} \,\, -\w_j \\
&=& - \max_{e \in C^*} \,\, \min_{f \in C^*} \,\, \w_f,
\end{eqnarray*}
remembering that the cocircuits of $M$ are precisely the
complements of its hyperplanes. This last expression is the result
of applying the blue rule to element $e$.
\end{proof}

\section{Phylogenetic trees} \label{section:trees}

Theorem \ref{main} is of particular interest when applied to
$M(K_n)$, the graphical matroid of the complete graph $K_n$. As
shown in \cite{Ardila}, the Bergman fan $\B(K_n)$ can be regarded
as a space of phylogenetic trees. Our results thus provide a new
point of view on a known algorithm in phylogenetics. In fact, it
is this context that provided the original motivation for our
results.

Let us now review the connection between the Bergman fan
$\B(K_n)$, ultrametrics, and phylogenetic trees. For more
information, see \cite{Ardila}.

\begin{definition}
A \emph{dissimilarity map} is a map $\delta: [n] \times [n]
\rightarrow \R$ such that $\delta(i,i)=0$ for all $i \in [n]$, and
$\delta(i,j) = \delta(j,i)$ for all $i,j \in [n]$. An
\emph{ultrametric} is a dissimilarity map such that, for all
$i,j,k \in [n]$, two of the values $\delta(i,j), \delta(j,k)$ and
$\delta(i,k)$ are equal and not less than the third.
\end{definition}

We can think of a dissimilarity map $\delta: [n] \times [n]
\rightarrow \R$ as a weight function $\w_{\delta} \in \R^{[n]
\choose 2}$ on the edges of the complete graph $K_n$. The
connection with our study is given by the following result.

\begin{theorem}\cite{Ardila}
A dissimilarity map $\delta$ is an ultrametric if and only if
$\w_{\delta}$ is a $M(K_n)$-ultrametric.
\end{theorem}

As mentioned earlier, the previous theorem is our reason for
giving vectors in $\B(M)$ the name of $M$-ultrametrics. We will
slightly abuse notation and write $\delta$ instead of
$\w_{\delta}$; this should cause no confusion.

\medskip

Ultrametrics are also in correspondence with a certain kind of
phylogenetic tree. Let $T$ be a rooted metric $n$-tree; that is, a
tree with $n$ leaves labelled $1,2,\ldots,n$, together with a
length assigned to each one of its edges. For each pair of leaves
$u,v$ of the tree, we define the \emph{distance} $d_T(u,v)$ to be
the length of the unique path joining leaves $u$ and $v$ in $T$.
This gives us a distance function $d_T:[n] \times [n] \rightarrow
\R$. We are interested in \emph{equidistant $n$-trees}. These are
the rooted metric $n$-trees such that the leaves are equidistant
from the root, and the lengths of the interior edges are positive.
(For technical reasons, the edges incident to a leaf are allowed
to have negative lengths.)

\begin{figure}[h]
\centering
\includegraphics[height=3cm]{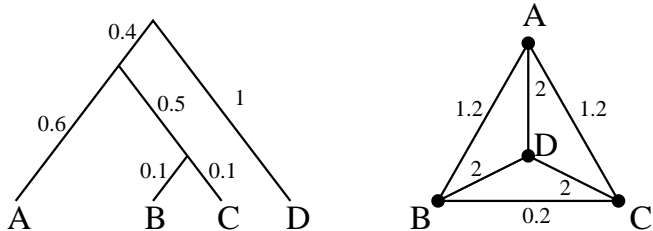}
\caption{An equidistant tree and its distance function.}
\label{fig:treeandultrametric}
\end{figure}

Figure \ref{fig:treeandultrametric} shows an example of an
equidistant $4$-tree, where the distance from each leaf to the
root is equal to $2$. It also shows the corresponding distance
function, recorded on the edges of the graph $K_4$. This distance
function is precisely the $M(K_4)$-ultrametric of Figure
\ref{fig:ultrametric}. This is not a coincidence, as the following
theorem shows.

\begin{theorem} \cite[Theorem 7.2.5]{Semple}
A dissimilarity map $\delta: [n] \times [n] \rightarrow \R$ is an 
ultrametric if and only if it is the distance function of an equidistant
$n$-tree.
\end{theorem}

We can think of equidistant trees as a model for the evolutionary
relationships between a certain set of species. The various
species, represented by the leaves, descend from a single root.
The descent from the root to a leaf tells us the history of how a
particular species branched off from the others, until the present
day. For more information on the applications of this and other
similar models, see for example \cite{Billera} and \cite{Semple}.

\medskip

One important problem in phylogenetics is the following: suppose
we have a way of estimating, in the present day, the pairwise
distances $d_T(i,j)$ between two species. The goal is to recover
the most likely tree $T$. It is well-known how to recover a tree
$T$ from its corresponding ultrametric $d_T$ \cite[Theorem
7.2.8]{Semple}. Of course, in real life, the measured data
$\delta(i,j)$ will not be exact. No tree will match the measured
distances exactly, and we need to find the tree which approximates
them most accurately. When proximity is measured in the
$\ell_{\infty}$ metric:
\[
||\delta-d_T||_{\infty} = \max |\,\delta(i,j)-d_T(i,j)|,
\]
Chepoi and Fichet gave a very nice answer to this question, which
we now review.

\medskip

Given a weight function $\w$ on the edges of $K_n$, let
\[
\w_U(x,y) = \min_{\textrm{paths $P$ from $x$ to $y$}} \,\,\,
\max_{\textrm{edges $e$ of $P$}} \w(e).
\]
It is not difficult to see that $\w_U$ is an ultrametric, known as
the \emph{subdominant ultrametric} of $\w$. Notice that $\w_U =
\w^{M(K_n)}$: the formula above is the result of applying the red
rule to edge $xy$.

Write $2\epsilon = ||\w - \w_U||_{\infty} =
\min{|\w(e)-\w_U(e)|}$, and define a second ultrametric by
$\w_U^{+\epsilon}(e) = \w_U(e) + \epsilon$ for each edge $e$ of
$K_n$.

\begin{theorem}\cite{Chepoi}
Given a dissimilarity map $\w$ on $[n]$, an
$\ell_{\infty}$-optimal ultrametric for $\w$ is the ultrametric
$\w_U^{+\epsilon}$.
\end{theorem}

Our discussion of matroid ultrametrics and tropical projection
provides a conceptual explanation of the previous theorem. Chepoi
and Fichet's goal is to construct an $\ell_{\infty}$-closest point
to $\w$ in the Bergman fan $\B(K_n)$. As mentioned in Section
\ref{section:tropical}, we can, in fact, essentially solve this
problem for any tropical polytope via tropical projection.

In tropical projective space, an $\ell_{\infty}$-optimal
ultrametric is $\pi_{\B(K_n)}(\w)$, the tropical projection of
$\w$ onto $\B(K_n)$. By Proposition \ref{projection}, this is also
the subdominant $M(K_n)$-ultrametric, $\w^{M(K_n)}$. Therefore, in
real space, an $\ell_{\infty}$-optimal ultrametric will be
$\w^{M(K_n)}$, up to a shift by a constant. It is easy to see that
the optimal constant is $\frac12 ||\w - \w^{M(K_n)}||_{\infty} =
\epsilon$.

\section{Acknowledgments}

I would like to thank the organizers and participants of the
Seminar on Phylogenetic Trees at the University of California,
Berkeley during the Fall of 2003, for a wonderful introduction to
this interesting subject. I would also like to thank Bernd
Sturmfels, who suggested the connection between tropical geometry
and Chepoi and Fichet's result, and Lior Pachter and Carly Klivans
for very helpful discussions.

\small{

}

\end{document}